\documentclass[12pt]{amsart}
\usepackage{amsmath, latexsym, amssymb, amsfonts,amscd,enumerate}
\usepackage{epstopdf}
\usepackage{graphicx}
\usepackage[myheadings]{fullpage}
\usepackage{epstopdf}
\usepackage{subfigure}
\usepackage[numbers]{natbib}
\theoremstyle{plain} \newtheorem{thm}{Theorem}[section]

\newtheorem{claim}[thm]{Claim}

\title{Two New Interpretations of the Fishburn Numbers and their Refined Generating Functions}
\author{Paul Levande}
\begin{document}
\maketitle
\begin{abstract}
We show that two classes of combinatorial objects--inversion tables with no subsequence of decreasing consecutive numbers and matchings with no 2-nestings--are enumerated by the Fishburn numbers.  In particular, we give a simple bijection between matchings with no 2-nestings and inversion tables with no subsequence of decreasing consecutive numbers.  We then prove using the involution principle that inversion tables with no subsequence of decreasing consecutive numbers have the same generating function as the Fishburn numbers.  The Fishburn numbers have previously been shown by Bousquet-M\'elou, Claesson, Dukes and Kitaev to enumerate $\textbf{(2+2)}$-avoiding posets, matchings with no left- or right-nestings, permutations avoiding a particular pattern, and so-called ascent sequences, and by Dukes and Parviainen to enumerate upper triangular matrices with non-negative entries and no empty rows or columns.  Claesson and Linusson conjectured they also enumerated matchings with no 2-nestings.  Using these new interpretations of the Fishburn numbers and another version of the involution, we prove the conjectured equality (also proven using matrices by Jel\'{\i}nek and by Yan) of two refinements by Remmel and Kitaev of the Fishburn generating function.  In an appendix, we state and prove another conjecture of Claesson and Linusson giving the distribution of left-nestings over the set of all matchings. 
\end{abstract}

\pagestyle{myheadings}
\markboth{Paul Levande}{New Interpretations of the Fishburn Numbers}
\section{Introduction}
It is well-known that the famous Catalan numbers $C_{n}$ enumerate, among many other classes of combinatorial objects, non-crossing and non-nesting matchings on $[2n]$, where a \textit{matching} is an involution with no fixed points, or a partition of $[2n]$ into disjoint pairs, and a \textit{nesting} in a matching $X$ is a pair of pairs $(a, b), (c, d) \in X$ such that $a < c < d < b$.  A superset of the set of non-nesting matchings on $[2n]$ is the set of non-neighbor-nesting matchings on $[2n]$, where a \textit{neighbor nesting} is a nesting $(a, b), (c, d) \in X$ such that either $c = a+1$ or $d = b-1$.  Zagier \cite{zaggen}, following Stoimenow \cite{stoen}, showed:
\begin{eqnarray}
\label{fishburndef}
\sum_{n=0}^{\infty} f_{n}t^{n} = 1+\sum_{m=1}^{\infty} \prod_{i=1}^{m} (1-(1-t)^{i})
\end{eqnarray}
where $f_{n}$ is the number of non-neighbor-nesting matchings on $[2n]$.  Recently, Bousquet-M\'elou et al \cite{BM}. showed that $f_{n}$ also enumerates other seemingly-disparate sets, each of which can be seen to be a superset of a set enumerated by the Catalan numbers:
\begin{itemize}
\item \textbf{(2+2)}-avoiding posets with $n$ elements, a superset of the set of \textbf{(2+2)}- and \textbf{(3+1)}-avoiding posets with $n$ elements enumerated by the $n$-th Catalan numbers $C_{n}$,
\item Permutations $\pi$ of $[n]$ such that, if $\pi = \pi_{1} \pi_{2} \ldots \pi_{n}$, there are no $i < j$ such that $\pi_{j} = \pi_{i}-1$ and $\pi_{i+1} > \pi_{i}$, a superset of the set of $231$-avoiding permutations of $[n]$ enumerated by the $n$-th Catalan number $C_{n}$.  
\item ascent sequences of length $n$, where an ascent sequence is a  sequence $x_{1} x_{2} \ldots x_{n}$ such that $x_{1} = 0$ and $0 \leq x_{i} \leq asc(x_{1}x_{2}\ldots x_{i-1})+1$, where $asc(x_{1} x_{2} \ldots x_{i-1})$ is the number of ascents of the sequence $x_{1} x_{2} \ldots x_{i-1}$, a superset of the set of sequences $x_{1} \ldots x_{n}$ such that $x_{1}=0$ and $x_{i+1} \leq x_{i}+1$ enumerated by the $n$-th Catalan number $C_{n}$.  
\end{itemize}
Claesson and Linusson \cite{claessonlinusson} recently conjectured another interpretation of the sequence $f_{n}$, which they refer to as the \textit{Fishburn numbers}: Matchings on $[2n]$ with no 2-nestings, where a $k$-nesting is a nesting $(a, b), (c, d)$ such that $a < c < d < b$ and $c-a \leq k$.  

Remmel and Kitaev \cite{remmellkit} recently gave the following refinement of (\ref{fishburndef}):
\begin{eqnarray}
\label{refineone}
\sum_{n=0}^{\infty}\sum_{d=1}^{n} f_{n, d}t^{n}z^{d} = 1+\sum_{m=1}^{\infty} \frac{zt}{(1-zt)^{m}}\prod_{i=1}^{m-1}(1-(1-t)^{i}),
\end{eqnarray}
where $f_{n, d}$ is (among other interpretations) the number of ascent sequences of length $n$ with $d$ zeroes.  Remmel and Kitaev also conjectured that the bivariate generating function $
\sum_{n=0}^{\infty}\sum_{d=1}^{n} f_{n, d}t^{n}z^{d}$ has the following simpler form:
\begin{eqnarray}
\label{refinetwo}
\sum_{n=0}^{\infty}\sum_{d=1}^{n} f_{n, d}t^{n}z^{d} = 1+\sum_{m=1}^{\infty} \prod_{i=1}^{m}(1-(1-t)^{i-1}(1-zt)).
\end{eqnarray}

In the following paper, we solve both conjectures and show how they are related: First, we give a simple bijection between matchings on $[2n]$ with no 2-nestings and inversion tables of length $n$ with no decreasing subsequence of consecutive numbers, i.e., sequences $a_{1} a_{2} \ldots a_{n}$ such that $0 \leq a_{i} \leq i-1$ for all $i$ and such that there exist no $p < q$ with $a_{p}=j+1$, $a_{q} = j$.  We then prove using the involution principle that the generating function of inversion tables of length $n$ with no decreasing subsequence of consecutive numbers is (\ref{fishburndef}).  We then show using the same proof that, if $f_{n, d}$ is now taken to be the number of such inversion tables  $a_{1} a_{2} \ldots a_{n}$ where $a_{i} = i-1$ for precisely $d$ distinct values of $i$, the bivariate generating function of $f_{n, d}$ is given by (\ref{refinetwo}).  Finally, we show using a variation of this proof that the the bivariate generating function of $f_{n, d}$, under this interpretation, is also given by (\ref{refineone}).  In an appendix, we give a brief solution to a related conjecture of Claesson and Linusson, that the distribution of \textit{left-nestings}, or nestings $(a, b), (c, d)$ with $c = a+1$, over all matchings is given by the second-order Eulerian triangle.

Note: Distinct proofs of the equality of (\ref{refineone}) and (\ref{refinetwo}), both using a matrix-based interpretation of the Fishburn numbers due to Dukes and Parviainen \cite{dukes}, were recently given independently by Jel\'{\i}nek \cite{Jel} and by Yan \cite{Yanproof}.  This work was done independently from the work presented here.  
\section{From factorial matchings to inversion tables}
\begin{figure}
\includegraphics[width=9 cm]{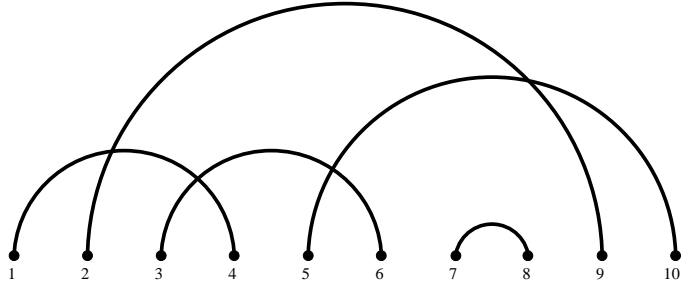}
\caption{A matching on $[10]$.}
\label{fishburnexam1}
\end{figure}

First, a note on visual notation: A matching $X$ on $[2n]$ can be illustrated by a diagram of $n$ semicircular arcs, where an arc connects $a$ and $b$ if and only if $(a, b) \in X$.  For example, the matching $(1, 4) (2, 9) (3, 6) (5, 10) (7, 8)$ on $[10]$ is illustrated by the diagram of $5$ semicircular arcs in Figure \ref{fishburnexam1}.
Note that an arc diagram is equivalent to its matching, i.e., the arc diagram does not have to be labelled to specify a unique matching.  
Claesson and Linusson \cite{claessonlinusson} proved the following claim about \textit{left-nestings}, where a left-nesting is a nesting $(a, b), (c, d)$ with $c = a+1$.  (For example, $(2, 9) (3, 6)$ in Figure \ref{fishburnexam1} is a left-nesting.)

\begin{claim} \label{inductclaim} The number of matchings on $[2n]$ with no left-nestings is $n!$. 
\end{claim}
\begin{proof} This is perhaps easiest to see visually and inductively: Let $X$ be a matching on $[2n]$ with no left-nestings.  A matching $X'$ on $[2n+2]$ can be formed from $X$ by adding an arc whose right endpoint is the rightmost part of the diagram, i.e., by adding $(a, 2n+2)$ to $X$, and re-labeling (moving all endpoints greater than $a$ over by one).  If the left endpoint is placed immediately before the left endpoint of another arc, i.e., if $(a, b) \in X$ and $a < b$, this will form a left-nesting.  Therefore, $X'$ will no left-nestings if and only if $(a, b) \in X$ and $a > b$, i.e., if and only if the left endpoint of the newest arc is be placed immediately before the right endpoint of an arc (including its own).  Since there will be $n+1$ right endpoints, there are $n+1$ possible matchings on $[2n+2]$ with no left-nestings formed by adding an arc to $X$.

Conversely, removing the pair $(a, 2n+2)$ from some matching $Z'$ on $[2n+2]$ with no left-nestings and re-labeling will produce a matching $Z$ with no left-nestings: The only left-nesting this could create in $Z$ is $(a-1, b), (a, c)$ with $a < c < b$, but then $(a+1, c+1), (a, 2n+2) \in Z'$ would be a left-nesting.  Therefore, every matching on $[2n+2]$ with no left-nestings is formed from adding an arc in $n+1$ possible ways to a matching on $[2n]$, and if there are $n!$ matchings on $[2n]$ with no left-nestings, there will be $(n+1)!$ matchings on $[2n+2]$ with no left-nestings.
\end{proof}

\begin{figure}
\includegraphics[width = 13 cm]{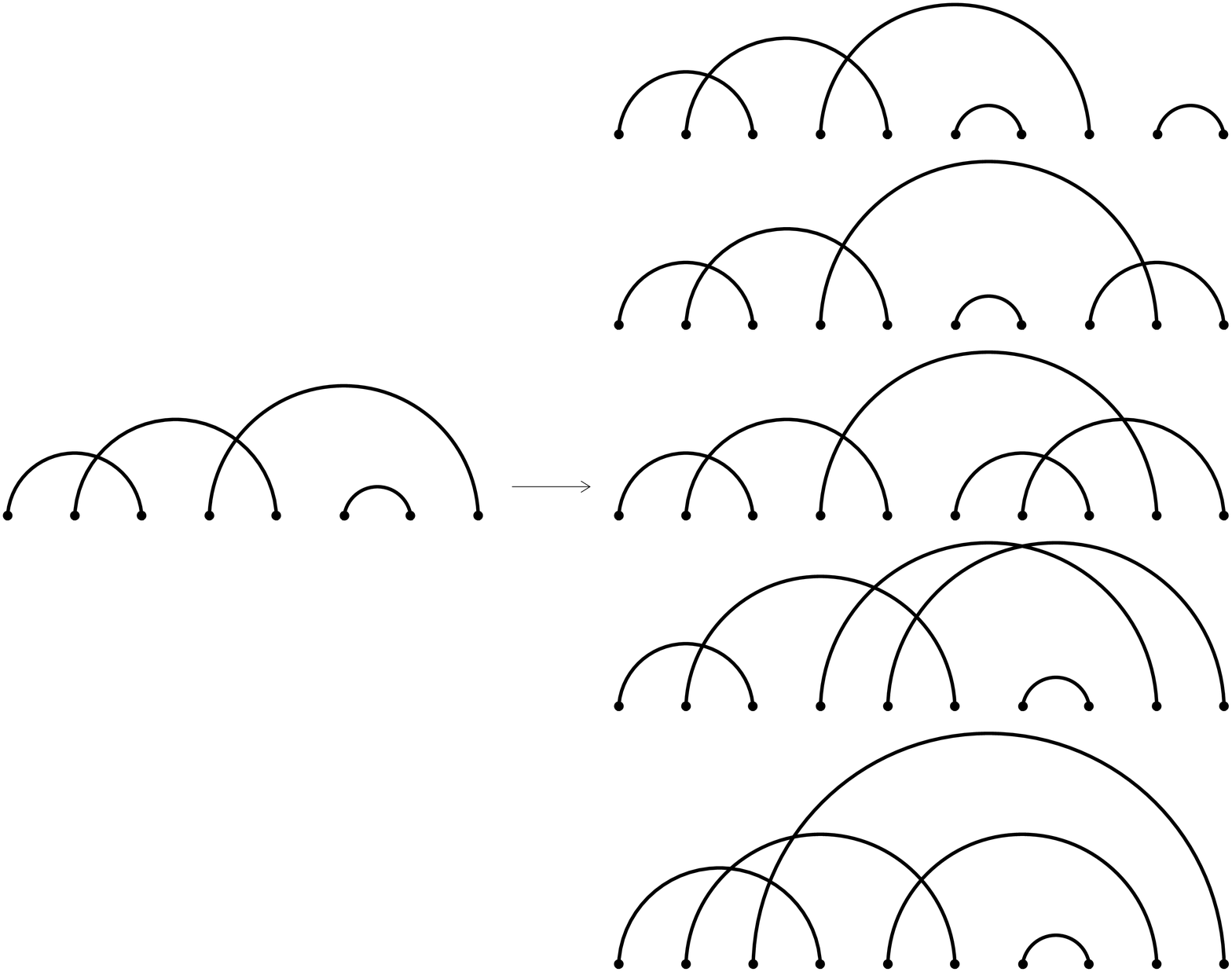}
\caption{Left-nesting-avoiding matchings on $[8]$ and $[10]$}
\label{fishburnexam2}
\end{figure}

We illustrate the inductive argument in Figure \ref{fishburnexam2}.

Note that we can (following Claesson and Linusson) refine the above proof: Let $X$ be a matching with no left-nestings, and let $\phi(X) = a_{1} a_{2} \ldots a_{n}$, where $a_{i}$ is the number of right endpoints to the left of the left endpoint of the $i$-th arc, where the arcs are numbered from left to right by right endpoint.  For example, if $X$ is the matching on $[8]$ in Figure \ref{fishburnexam2}, $\phi(X) = 0021$.
\begin{claim}
$\phi$ gives a bijection between matchings with no left-nestings and inversion tables, i.e., sequences $a_{1} a_{2} \ldots a_{n}$ such that $0 \leq a_{i} \leq i-1$ for all $1 \leq i \leq n$.
\end{claim}
\begin{proof} This should be clear from the above: A new arc added to a matching on $[2n]$ corresponding to a pair $(a, 2n+2)$ can have its left endpoint be placed only immediately before the first right endpoint, the second right endpoint, and so on, i.e., to the right of $0, 1, \ldots, n$ right endpoints.
\end{proof}
For example, when applied to the the matchings on $[10]$ in Figure \ref{fishburnexam2}, $\phi$ gives, from top to bottom, $00214, 00213, 00212, 00211$, and $00210$.

In order to prove that the number of matchings on $[2n]$ with no $2$-nestings is the $n$-th Fishburn number $f_{n}$, we restrict $\phi$ to matchings with no $2$-nestings:

\begin{claim} $\phi$ gives a bijection between matchings with no $2$-nestings and inversion tables $a_{1} a_{2} \ldots a_{n}$ with no decreasing subsequences of consecutive integers, i.e., with no $p < q$ such that $a_{p} = j+1$ and $a_{q} = j$ for some $j$. 
\end{claim}
\begin{proof} Recall that a $2$-nesting is a nesting $(a, b), (c, d)$ with $c-a \leq 2$.  In particular, if $X$ is a matching with no left-nestings and $\phi(X) = a_{1} a_{2} \ldots a_{n}$, then $X$ will have a $2$-nesting if and only if, for some $a$, $(a, b), (a+2, d) \in X$, $d < b$, and, if $(a+1, g) \in X$, $g < a$: if $a+2 < g < b$, then $(a, b)$ and $(a+1, g)$ will be a left-nesting, and if $b < g$, $(a+1, g), (a+2, d)$ will be a left-nesting.  Therefore $X$ will have a $2$-nesting if and only if its arc diagram gives precisely one right endpoint between the left endpoints of the arc corresponding to $(a, b)$ and the arc corresponding to $(a+2, d)$, or if and only if $a_{p} = j+1$ and $a_{q} = j$, where $d$ is the $p$-th right endpoint, $b$ is the $q$-th right endpoint, and $j$ is the number of right endpoints to the left of $a$. 

Similarly, if, for some $p < q$, $a_{p} = j+1$ and $a_{q} = j$, then the arc diagram of $X$ will have precisely one right endpoint in between the left endpoint of the $q$-th arc and the left endpoint of the $p$-th arc.  Let the $q$-th arc correspond to the pair $(e, f) \in X$ and the $p$-th arc correspond to the pair $(g, h) \in X$, with the single arc with right endpoint in between correspond to the pair $(r, s) \in X$, so $e < s < g < h < f$.  $e+1, e+2, \ldots, s-1, s+2, s+3, \ldots, g-1$ must all be left endpoints of arcs, since there is only one right endpoint in between $e$ and $g$.   In general, let $X(a)$ be the element paired with $a$ in $X$.  
To avoid left-nestings, $X(e+1) > X(e) = f$, $X(e+2) > X(e+1)$, and so on, so $X(s-1) > X(s-2) > \ldots > X(e) = f$  Similarly, to avoid left-nestings, $X(g-1) < X(g) = h$, $X(g-2) < X(g-1)$, and so on, so $X(s+1) < X(s+2) < \ldots < X(g) =h$.  Therefore $(s-1, X(s-1)), (s+1, X(s+1)) \in X$ and $X(s+1) < h < f < X(s-1)$, and so $(s-1, X(s-1)), (s+1, X(s+1)) \in X$ is a 2-nesting, and $X$ has at least one $2$-nesting if and only if $\phi(X)$ has at least one decreasing subsequence of consecutive integers.  
\end{proof}
Let $T_{n}$ be the set of inversion tables of length $n$ with no decreasing subsequences of consecutive integers.  
\section{The main proof}
We will now prove that $|T_{n}| = f_{n}$, i.e., that the set of inversion tables of length $n$ with no decreasing subsequences of consecutive integers is enumerated by the $n$-th Fishburn number.  In particular, we will prove that:
\begin{eqnarray}
\label{fishburntable}
\sum_{n=0}^{\infty} |T_{n}|t^{n} = 1+\sum_{m=1}^{\infty} \prod_{i=1}^{m} (1-(1-t)^{i}).
\end{eqnarray}
First, we will define a class of diagrams to give a visual interpretation to the right-hand-side of (\ref{fishburntable}).  Given a staircase partition diagram $(1, 2, \ldots, k)$, a \textit{filling} of the diagram is a placement of dots into the squares of the diagram, with at most one dot per square.  We say that a diagram with filling is a \textit{Fishburn diagram} if every column has at least one dot.  Let $Y$ be the set of Fishburn diagrams of any size, with $Y_{n}$ the set of Fishburn diagrams with $n$ dots.  Equivalently, let $Y$ be the set of sequences of sets $A_{1} A_{2} \ldots A_{k}$, where $k$ can vary, such that $A_{i} \subset \left\{0, 1, \ldots, i-1 \right\}$ and $A_{i} \neq \emptyset$, with $Y_{n}$ the set of such sequences such that $\sum_{i=1}^{k} |A_{i}| = n$.  We will refer to Fishburn diagrams and sequences of sets interchangably; for the correspondence, given a sequence of sets $A_{1} A_{2} \ldots A_{k}$ with the above condition, the corresponding diagram is a filling of $(1, 2, \ldots, k)$ with a dot in the $j$-th row from the bottom in the $i$-th column from the left if and only if $j-1 \in A_{i}$.  

\begin{claim}
The following identity of generating functions holds:
\begin{eqnarray}
\label{fishburnweight}
\sum_{n=0}^{\infty} \sum_{A_{1} A_{2} \ldots A_{k} \in Y_{n}} t^{n}(-1)^{n-k} =  1+\sum_{m=1}^{\infty} \prod_{i=1}^{m} (1-(1-t)^{i}).
\end{eqnarray}
\end{claim}
\begin{proof}
Let us first consider the weighted sum of Fishburn diagrams of length $m$, i.e., sequences $A_{1} A_{2} \ldots A_{m} \in Y$, beginning by giving each column a weight of $-1$.  For a specific column $A_{i}$ with $1 \leq i \leq m$, we can either place a dot, or not, each square.  Let the weight of a dot be $-t$ and the weight of a square with no dot be $1$.  Then for each square, our weighted choice is $1-t$, and there are $i$ squares in this column.  Since we must pick at least one dot for the column, the weighted sum over all possible choices is $1-(1-t)^{i}$.
Therefore,
\begin{eqnarray*}
\sum_{A_{1}A_{2} \ldots A_{m} \in Y} (-t)^{|A_{1}|+|A_{2}|+\ldots+|A_{m}|} (-1)^{m} = \sum_{m=1}^{\infty}\prod_{i=1}^{m} (1-(1-t)^{i}),
\end{eqnarray*}
and
\begin{eqnarray}
\label{fishburnlem}
\sum_{m=0}^{\infty}\sum_{A_{1}A_{2} \ldots A_{m} \in Y} (-t)^{|A_{1}|+|A_{2}|+\ldots+|A_{m}|} (-1)^{m} = 1+\sum_{m=1}^{\infty}\prod_{i=1}^{m} (1-(1-t)^{i}).
\end{eqnarray} 
It should be clear that (\ref{fishburnlem}) is equivalent to (\ref{fishburnweight}).
\end{proof}

\begin{figure}
\label{fishburnexam3}
\includegraphics{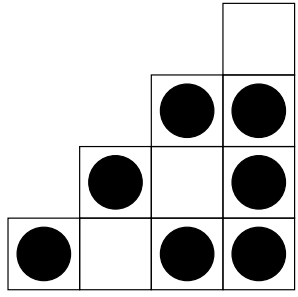}
\caption{A Fishburn diagram in $Y_{7}$ with signed weight $(-1)^{3}t^{7}$.}
\end{figure}

We can think of the left-hand-side of (\ref{fishburnweight}) as the signed weighted sum over all Fishburn diagrams of any size, where the weight is given by the number of dots and the sign is given by the number of ``extra'' dots, i.e., dots other than the minimal one per column.  For example, Figure \ref{fishburnexam3} shows the Fishburn diagram corresponding to $A_{1} A_{2} A_{3} A_{4} \in Y_{7}$, where $A_{1} = \left\{0 \right\}, A_{2} = \left\{1 \right\}, A_{3} = \left\{0, 2 \right\}$ and $A_{4} = \left\{0, 1, 2\right\}$, which will have a signed weight of $(-1)^{3}t^{7}$.

In order to prove (\ref{fishburntable}), we will define an involution $\psi_{n}: Y_{n} \rightarrow Y_{n}$ such that:
\begin{itemize}
\item $A_{1} A_{2} \ldots A_{k} \in Fix(\psi_{n})$ if and only if, for all $1 \leq i \leq k$, $|A_{i}| = 1$ (so $k$ = $n$) and there are no $p, q, j$ with $p < q$, $A_{p} = \left\{j+1 \right\}, A_{q} = \left\{j \right\}$,
\item If $A_{1} A_{2} \ldots A_{k} \notin Fix(\psi_{n})$ and $\psi_{n}(A_{1} A_{2} \ldots A_{k}) = B_{1} B_{2} \ldots B_{r}$, $r = k \pm 1$,
\end{itemize}
where $Fix(\psi_{n}) = \left\{A_{1} A_{2} \ldots A_{k} \in Y_{n} : \psi_{n}(A_{1} A_{2} \ldots A_{k}) = A_{1} A_{2} \ldots A_{k} \right\}$.  If we can define such an involution $\psi_{n}$, then we will have
\begin{eqnarray*}
\sum_{n=0}^{\infty} |Fix(\psi_{n})|t^{n}= 1+\sum_{m=1}^{\infty}\prod_{i=1}^{m} (1-(1-t)^{i}).
\end{eqnarray*}
This will prove (\ref{fishburntable}), since $T_{n}$ will be in trivial bijection with $Fix(\psi_{n})$.  

Given $A_{1} A_{2} \ldots A_{k} \in Y_{n}$, let $j$ be the smallest integer such that either:
\begin{itemize}
\item For some $i$, $j \in A_{i}$ and $|A_{i}| > 1$, or
\item For some $p < q$, $j+1 \in A_{p}$, $j \in A_{q}$.  
\end{itemize}
If there is no such $j$, define $\psi_{n}(A_{1} A_{2} \ldots A_{k}) = A_{1} A_{2} \ldots A_{k}$.  Otherwise, we divide into cases:
\begin{description}
\item[Case 1] If, for some $i$, $j \in A_{i}$ and $|A_{i}| > 1$, let $I$ be the minimal such $i$.  By the minimality of $j$, $j$ is the smallest member of $A_{I}$.  Let $j+R$ be the second-smallest member of $A_{I}$.  Define sets $B_{L}$ as follows for $1 \leq L \leq k+1$.  
\begin{displaymath}
B_{L} = \begin{cases}
A_{L} & 1 \leq L \leq I-R \\
\left\{s-R+1 : s \in A_{I}, s \neq j \right\} & L = I-R+1 \\
\left\{s: s \in A_{L-1}, s < j+1 \right\} \bigcup \left\{s+1: s \in A_{L-1}, s \geq j+1 \right\} & I-R+2 \leq L \leq I \\
\left\{j \right\} & L = I+1 \\
\left\{s: s \in A_{L-1}, s < j \right\} \bigcup \left\{s+1: s \in A_{L-1}, s \geq j \right\}& I+2  \leq L \leq k+1 
\end{cases}
\end{displaymath}
It should be clear that $(B_{1}, B_{2}, \ldots, B_{k+1}) \in Y_{n}$.  Define $\psi_{n}(A_{1}, A_{2}, \ldots, A_{k}) =   (B_{1}, B_{2}, \ldots, B_{k+1})$.

Note that, by construction, $B_{I+1} = \left\{j \right\}$ is the last set of $B_{1}, B_{2}, \ldots, B_{k+1}$ to contain $j$, and $B_{I-R+1}$ is the last set of $B_{1}, B_{2}, \ldots, B_{I+1}$ to contain $j+1$, which must be the smallest member of $B_{I-R+1}$.  Note also that, since $A_{I}$ was the first set of $A_{1}, A_{2}, \ldots, A_{k}$ to include $j$ and have multiple members, no member of $B_{1}, B_{2}, \ldots, B_{k+1}$ will contain $j$ and have multiple members.

\item[Case 2] If there is no $i$ such that $j \in A_{i}$ and $|A_{i}| > 1$, then, for some $p < q$, $j+1 \in A_{p}$ and $j \in A_{q}$.  Let $Q$ be maximal such that $j \in A_{Q}$.  By definition, $A_{Q} = \left\{j \right\}$.  Let $P$ be maximal such that $j+1 \in A_{P}$ and $P < Q$.  Therefore $j+1 \notin A_{L}$ for $P+1 \leq L \leq Q-1$ and $j \notin A_{L}$ for $Q+1 \leq L \leq k$.  By the minimality of $j$, $j+1$ must be the smallest member of $A_{P}$.  Define sets $B_{L}$ as follows for $1 \leq L \leq k-1$:
\begin{displaymath}
B_{L} = 
\begin{cases}
A_{L} & 1 \leq L \leq P-1 \\
\left\{s: s \in A_{L+1}, s < j+1 \right\} \bigcup \left\{s-1 : s \in A_{L+1}, s > j+1 \right\} & P \leq L \leq Q-2 \\
\left\{s+M-N-1: s \in A_{P} \right\} \bigcup \left\{j \right\} & L = Q-1 \\
\left\{s: s \in A_{L+1}, s < j \right\} \bigcup \left\{s-1: s \in A_{L+1}, s > j \right\} & Q \leq L \leq k-1
\end{cases}
\end{displaymath}
It should be clear that $(B_{1}, B_{2}, \ldots, B_{k-1}) \in Y_{n}$.  Define $\psi_{n}(A_{1}, A_{2}, \ldots, A_{k}) =   (B_{1}, B_{2}, \ldots, B_{k-1})$.  

Note that, by construction, $B_{Q-1}$ will be the first set of $B_{1}, B_{2}, \ldots, B_{k-1}$ to contain $j$ and have multiple members, with $j+Q-P$ the second-smallest member.   
\end{description}
Note as well that, in either case, $\psi_{n}$ preserves the minimality of $j$, and that therefore $\psi_{n}$ is an involution.  This proves (\ref{fishburntable}).

\begin{figure}
\label{fishburnexam4}
\includegraphics[width = 7 cm]{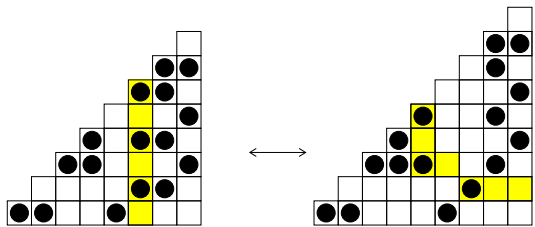}
\caption{$\psi_{11}: Y_{11} \rightarrow Y_{11}$}
\end{figure}

We illustrate the involution in Figure \ref{fishburnexam4}.  The Fishburn diagram on the left corresponds to $A_{1} A_{2} A_{3} A_{4} A_{5} A_{6} A_{7} A_{8} \in Y_{16}$, where $A_{1} = A_{2} = A_{5} = \left\{0 \right\}$, $A_{3} = \left\{2 \right\}$, $A_{4} = \left\{2, 3 \right\}$, $A_{6} = \left\{1, 3, 5 \right\}$, $A_{7} = \left\{1, 3, 5, 6 \right\}$, and $A_{8} = \left\{2, 4, 6 \right\}$.  It will have a signed weight of $(-1)^{8}t^{16}$.  The column corresponding to $\left\{1, 3, 5 \right\}$ has been highlighted.  The Fishburn diagram on the right oorresponds to $B_{1} B_{2} B_{3} B_{4} B_{5} B_{6} B_{7} B_{8} B_{9} \in Y_{16}$, where $B_{1} = B_{2} = B_{6} = \left\{0 \right\}$, $B_{3} = \left\{2 \right\}$, $B_{4} = \left\{2, 3 \right\}$, $B_{5} = \left\{2, 4 \right\}$, $B_{7} = \left\{1 \right\}$, $B_{8} = \left\{2, 4, 6, 7 \right\}$, and $B_{9} = \left\{3, 5, 7 \right\}$.  It will have a signed weight of $(-1)^{7}t^{16}$.  Squares have been highlighted to indicate how $2, 4 \in B_{5}$, $2 \notin B_{6}$, $1 \in B_{7}$ and $1 \notin B_{8}, B_{9}$.  Informally, the dots corresponding to $3$ and $5$ in $A_{6}$ have been ``moved'' into dots corresponding to $2$ and $4$ in what is now the column $B_{5}$, leaving the dot corresponding to the $1$ in $A_{6}$ now in a column $B_{7}$ one square taller, with squares added in the third row from the bottom in $B_{6}$ and the second row from the bottom in $B_{8}$ and $B_{9}$. 

\section{Refinements of the generating function of the Fishburn numbers}
Let $T_{n, d}$ be the subset of $T_{n}$ consisting of inversion tables $a_{1} a_{2} \ldots a_{n} \in T_{n}$ such that $a_{i} = i-1$ for exactly $d$ distinct values of $i$.  We will prove the following two generating function identities,
\begin{eqnarray}
\label{refinetableone}
\sum_{n=0}^{\infty} \sum_{d=1}^{n} |T_{n, d}|t^{n}z^{d} = 1+\sum_{m=1}^{\infty} \prod_{i=1}^{m}(1-(1-t)^{i-1}(1-zt)),
\end{eqnarray}
and
\begin{eqnarray}
\label{refinetabletwo}
\sum_{n=0}^{\infty} \sum_{d=1}^{n} |T_{n, d}|t^{n}z^{d} = 1+\sum_{m=1}^{\infty} \frac{zt}{(1-zt)^{m}}\prod_{i=1}^{m-1}(1-(1-t)^{i}),
\end{eqnarray}
using variations on the proof of (\ref{fishburntable}).  This will prove the equality of (\ref{refineone}) and (\ref{refinetwo}).  
\subsection{The simpler refinement}
We begin with a proof of (\ref{refinetableone}).  Let $Y_{n, d}$ be the subset of $Y_{n}$ consisting of Fishburn diagrams $A_{1} A_{2} \ldots A_{k} \in Y_{n}$ such that $i-1 \in A_{i}$ for exactly $d$ distinct values of $i$.  
\begin{claim}
\begin{eqnarray}
\label{refineclaim}
\sum_{n=0}^{\infty} \sum_{d=1}^{n} \sum_{A_{1} A_{2} \ldots A_{k} \in Y_{n, d}} z^{d}t^{n}(-1)^{n-k} = 1+\sum_{m=1}^{\infty} \prod_{i=1}^{m}(1-(1-t)^{i-1}(1-zt)).
\end{eqnarray}
\end{claim}
\begin{proof}
We follow the proof of (\ref{fishburnweight}) and first consider the weighted sum of Fishburn diagrams of length $m$, i.e., sequences $A_{1} A_{2} \ldots A_{m} \in Y$, beginning by giving each column a weight of $-1$.  For a specific column $A_{i}$ with $1 \leq i \leq m$, we can either place a dot, or not, each square.  Let the weight of a dot below the $i$-th row from the bottom be $-t$, the weight of a dot in the $i$-th row from the bottom be $-zt$, and the weight of a square with no dot be $1$.  Then for each square, our weighted choice is $1-t$, except for the top square in the column, where the weighted choice is $1-zt$, and there are $i-1$  other squares in this column.  Since we must pick at least one dot for the column, the weighted sum over all possible choices is $1-(1-t)^{i-1}(1-zt)$, and it should be clear this proves (\ref{refineclaim}).
\end{proof}

\begin{claim}
$\psi_{n}$ restricts to an involution from $Y_{n, d} \rightarrow Y_{n, d}$, i.e., if $A_{1} A_{2} \ldots A_{k} \in Y_{n, d}$, then $\psi_{n}(A_{1} A_{2} \ldots A_{k}) \in Y_{n, d}$.
\end{claim}
\begin{proof} If $A_{1} A_{2} \ldots A_{k} \in Fix(\psi_{n})$, then this is clearly true.  Otherwise, let $j$ be minimal such that either some $j \in A_{i}$ and $|A_{i}| > 1$ or $j+1 \in A_{p}, j \in A_{q}$ for $p < q$.  Assume the former holds, i.e., $j \in A_{i}$ and $|A_{i}| > 1$.  Let $I$ be the minimal such $i$, $j+R \in A_{I}$ be the second-smallest element of $A_{I}$, and let $\psi(A_{1} A_{2} \ldots A_{k}) = B_{1} B_{2} \ldots B_{k+1}$.  From the definition of $\psi_{n}$, we see that $L-1 \in B_{L}$ if and only if
\begin{eqnarray*}
L-1 \in A_{L} & 1 \leq L \leq I-R \\
I-1 \in A_{I} &  L = I-R+1 \\
L-2 \in A_{L-1} & I-R+2 \leq L \leq k+1,  L \neq I+1
\end{eqnarray*}
Therefore $B_{1} B_{2} \ldots B_{k+1} \in Y_{n, d}$.  Since $\psi_{n}$ is an involution, this suffices.
\end{proof}.

For example, both of the Fishburn diagrams in Figure \ref{fishburnexam4} are in $Y_{11, 5}$.  Let $\psi_{n}|_{Y_{n, d}} = \psi_{n, d}$.  Then $\psi_{n, d}: Y_{n, d} \rightarrow Y_{n, d}$ is an involution, and therefore 
\begin{eqnarray}
\label{refineclaimtwo}
\sum_{n=0}^{\infty} \sum_{d=1}^{n} |Fix(\psi_{n, d})|z^{d}t^{n} = 1+\sum_{m=1}^{\infty} \prod_{i=1}^{m}(1-(1-t)^{i-1}(1-zt)).
\end{eqnarray}
It should be clear that (\ref{refineclaimtwo}) is equivalent to (\ref{refinetableone}), since $Fix(\psi_{n, d}) = Fix(\psi_{n}) \bigcap Y_{n, d}$ is in trivial bijection with $T_{n, d}$.  
\subsection{The more complicated refinement}
We will now prove (\ref{refinetabletwo}).  First, we will define a set $\widetilde{Y_{n, d}}$ which has the right-hand-side of (\ref{refinetabletwo}) as a natural generating function.  Second, we will show that the right-hand-side of (\ref{refinetabletwo}) is equal to the bivariate generating function of the fixed point set of an involution $\widetilde{\psi_{n, d}}: \widetilde{Y_{n, d}} \rightarrow \widetilde{Y_{n, d}}$.  Finally, we will define a natural bijection from $Fix(\widetilde{\psi_{n, d}})$ to $Fix(\psi_{n, d})$. 

Let $\widetilde{Y_{n, d}}$ be the set of ordered pairs $(\lambda, A_{1} A_{2} \ldots A_{k})$ where $\lambda = (\lambda_{1}, \lambda_{2}, \ldots, \lambda_{k+1})$ is a composition of $d$ into $k+1$ non-negative parts, with $\lambda_{1} > 0$, and $A_{1} A_{2} \ldots A_{k} \in Y_{n-d}$.  
\begin{claim}
\begin{eqnarray}
\label{refineclaimthree}
\sum_{n=0}^{\infty} \sum_{d=1}^{\infty}\sum_{(\lambda, A_{1} A_{2} \ldots A_{k}) \in \widetilde{Y_{n, d}}} z^{d}t^{n}(-1)^{n-k} = 1+\sum_{m=1}^{\infty} \frac{zt}{(1-zt)^{m}}\prod_{i=1}^{m-1}(1-(1-t)^{i}).
\end{eqnarray}
\end{claim}
\begin{proof}
Let us consider all $(\lambda, A_{1} A_{2} \ldots A_{m-1})$, where $m$ is fixed.  We weight $(\lambda, A_{1} A_{2} \ldots A_{m-1})$ by the product of the weight of $A_{1} A_{2} \ldots A_{m-1}$, as in the proof of (\ref{fishburnweight}) (but not the weight from the proof of (\ref{refineclaim})) and the weight of $\lambda$, which we take to be $(zt)^{|\lambda|}$, where $|\lambda| = \lambda_{1}+\lambda_{2}+\ldots+\lambda_{m}$.  The generating function of compositions with $m$ non-negative parts, where the first part must be positive, is $\frac{zt}{(1-zt)^{m}}$.  Then the weighted sum over all $(\lambda, A_{1} A_{2} \ldots A_{m-1})$ is $\frac{zt}{(1-zt)^{m}}\prod_{i=1}^{m-1}(1-(1-t)^{i})$, and (\ref{refineclaimthree}) follows.
\end{proof}
We will now define an involution $\widetilde{\psi_{n, d}}: \widetilde{Y_{n, d}} \rightarrow \widetilde{Y_{n, d}}$ such that
\begin{itemize}
\item $(\lambda, A_{1} A_{2} \ldots A_{k}) \in Fix(\psi_{n})$ if and only if, for all $1 \leq i \leq k$, $|A_{i}| = 1$ (so $k$ = $n-d$) and if for some $j$, $A_{p} = \left\{j+1 \right\}, A_{q} = \left\{j \right\}$, and $p < q$, then $\lambda_{j+2} \neq 0$,
\item If $(\lambda, A_{1} A_{2} \ldots A_{k}) \notin Fix(\widetilde{\psi_{n,d}})$ and $\widetilde{\psi_{n, d}}(\lambda, A_{1} A_{2} \ldots A_{k}) = (\mu, B_{1} B_{2} \ldots B_{r})$, $r = k \pm 1$.
\end{itemize}

Given $(\lambda, A_{1} A_{2} \ldots A_{k}) \in \widetilde{Y_{n, d}}$, let $j$ be the smallest integer such that either, for some $i$,  $j \in A_{i}$ and $|A_{i}| > 1$ or, for some $p, q$, $j+1 \in A_{p}, j \in A_{q}$, and $\lambda_{j+2} = 0$.  If there is no such $j$, define $\widetilde{\psi_{n, d}}(\lambda, A_{1} A_{2} \ldots A_{k}) = (\lambda, A_{1} A_{2} \ldots A_{k})$.  Otherwise, 
\begin{description}
\item[Case 1] If, for some $i$, $j \in A_{i}$ and $|A_{i}| > 1$, define $\widetilde{\psi_{n, d}}(\lambda, A_{1} A_{2} \ldots A_{k}) = (\mu, \psi_{n-d}(A_{1} A_{2} \ldots A_{k}))$, where $\mu = (\lambda_{1}, \lambda_{2}, \ldots, \lambda_{j+1}, 0, \lambda_{j+2}, \lambda_{j+3}, \ldots, \lambda_{k+1})$.  Note that $\psi_{n-d}(A_{1} A_{2} \ldots A_{k}) = B_{1} B_{2} \ldots B_{k+1} \in Y_{n-d}$ by the definition of $\psi_{n-d}$, so $(\mu, \psi_{n-d}(A_{1} A_{2} \ldots A_{k})) \in \widetilde{Y_{n, d}}$.  

\item[Case 2] If there is no $i$ such that $j \in A_{i}$ and $|A_{i}| > 1$, then, for some $p < q$, $j+1 \in A_{p}$, $j \in A_{q}$, and $\lambda_{j+2} = 0$.  Define $\widetilde{\psi_{n, d}}(\lambda, A_{1} A_{2} \ldots A_{k}) = (\mu, \psi_{n-d}(A_{1} A_{2} \ldots A_{k}))$, where $\mu = (\lambda_{1}, \lambda_{2}, \ldots, \lambda_{j+1}, \lambda_{j+3}, \ldots, \lambda_{k+1})$.  Note that $\psi_{n-d}(A_{1} A_{2} \ldots A_{k}) = B_{1} B_{2} \ldots B_{k-1} \in Y_{n-d}$ by the definition of $\psi_{n-d}$, so $(\mu, \psi_{n-d}(A_{1} A_{2}, \ldots A_{k})) \in \widetilde{Y_{n, d}}$.
\end{description}

Therefore, 
\begin{eqnarray} \label{refineclaimfour}
\sum_{n=0}^{\infty} \sum_{d=1}^{n} |Fix(\widetilde{\psi_{n, d}})|t^{n}z^{d} = 1+\sum_{m=1}^{\infty} \frac{zt}{(1-zt)^{m}}\prod_{i=1}^{m-1}(1-(1-t)^{i}).
\end{eqnarray}

To finish the proof of (\ref{refinetabletwo}), we need only define a bijection from $Fix(\widetilde{\psi_{n, d}})$ to $Fix(\psi_{n, d})$.

In fact, we claim that $\widetilde{Y_{n, d}}$ is in natural fixed-point-preserving bijection with a subset of $Y_{n, d}$.  Define an injection $f: \widetilde{Y_{n, d}} \rightarrow Y_{n, d}$ as follows: Given $(\lambda, A_{1} A_{2} \ldots A_{k}) \in \widetilde{Y_{n, d}}$, where $\lambda = (\lambda_{1}, \lambda_{2}, \ldots, \lambda_{k+1})$, begin with a staircase diagram $(1, 2, \ldots, d+k)$.  Place an $X$ in the top of the column $L$ if $L = \lambda_{1}+1, \lambda_{1}+\lambda_{2}+2, \lambda_{1}+\lambda_{2}+\lambda_{3}+3$, and so on up to $\lambda_{1}+\lambda_{2}+\ldots+\lambda_{k}+k$.  Place a dot in the top of all other columns; note that you have placed $d$ dots at the top of their columns, with $k$ other squares dividing the top diagonal into the composition $(\lambda_{1}, \lambda_{2}, \ldots, \lambda_{k+1})$.  Place an $X$ in each square in the same column as a dot at the top of its column, and in each square to the right of, and in the row immediately below, a dot at the top of its column.  The squares that now have neither a dot nor an $X$ will form a staircase sub-diagram with shape$(1, 2, \ldots, k)$.  Place the Fishburn diagram $A_{1} A_{2} \ldots A_{k}$, which has $n-d$ dots, into this staircase sub-diagram.  Let $D_{1} D_{2} \ldots D_{d+k}$ be the resulting Fishburn diagram, and define $f(\lambda, A_{1} A_{2} \ldots A_{k}) = D_{1} D_{2} \ldots D_{d+k}$.  

\begin{figure}
\includegraphics[width = 7 cm]{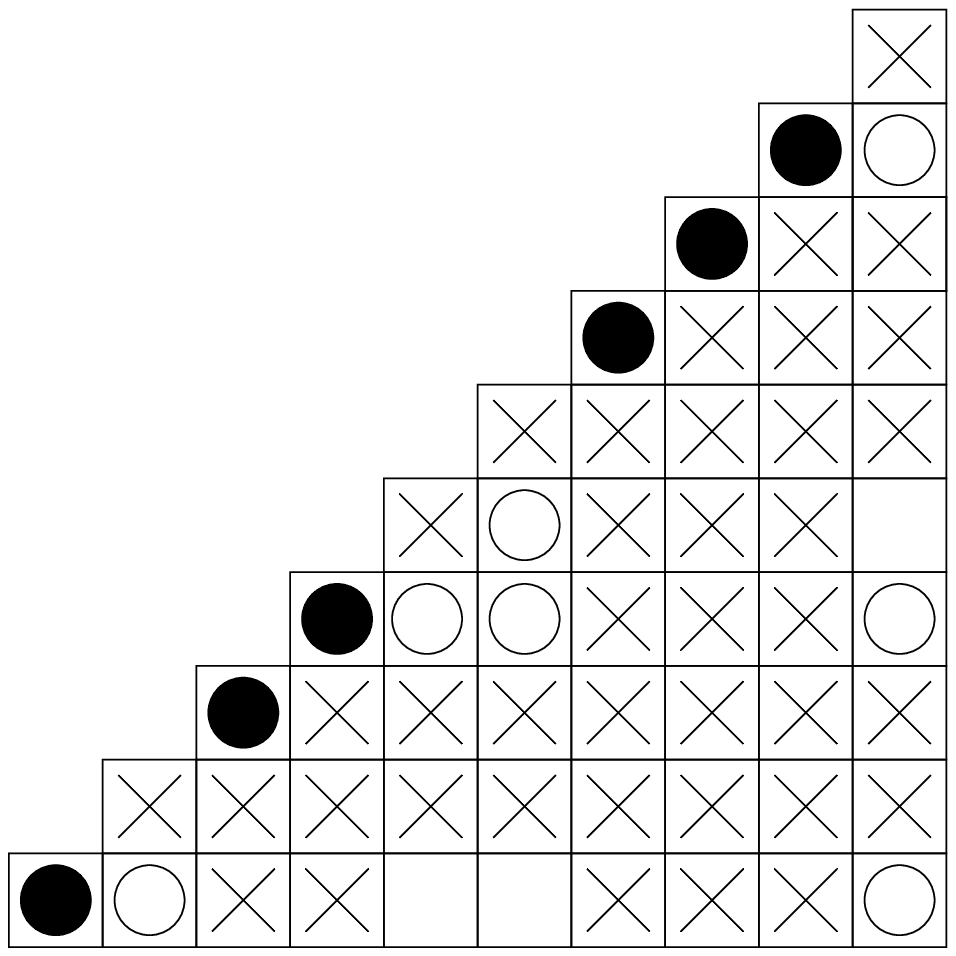}
\caption{A Fishburn diagram in $Y_{13, 6}$ corresponding to an element of $\widetilde{Y_{13, 6}}$.} 
\label{fishburnexam5}
\end{figure}

By construction, $D_{1} D_{2} \ldots D_{d+k}$ is a Fishburn diagram with $n-d+d = n$ dots, $d$ of which are at the tops of their columns, so $D_{1} D_{2} \ldots D_{d+k} \in Y_{n, d}$.  For example, if $\lambda = (1, 2, 0, 3, 0)$, $A_{1} = \left\{0 \right\}, A_{2} = \left\{1 \right\}, A_{3} = \left\{1, 2 \right\}$, and $A_{4} = \left\{0, 1, 3 \right\}$, then $\lambda$ is a composition of $6$ into $5$ parts and $A_{1} A_{2} A_{3} A_{4} \in Y_{7}$.  Therefore $(\lambda, A_{1} A_{2} A_{3} A_{4}) \in \widetilde{Y_{13, 6}}$, and $f(\lambda, A_{1} A_{2} A_{3} A_{4}) \in Y_{13, 6}$ is the Fishburn diagram in Figure \ref{fishburnexam5}.

(In Figure \ref{fishburnexam5} we used white for the seven dots from $A_{1}A_{2}A_{3}A_{4}$ and black for the six dots from $\lambda$.  Note how the $10$ squares with neither an $X$ nor a black dot form a subdiagram with staircase shape $(1, 2, 3, 4)$.). 

\begin{claim} $f$ is a bijection between $\widetilde{Y_{n, d}}$ and the subset of $Y_{n, d}$ consisting of Fishburn diagrams $D_{1} D_{2} \ldots D_{k}$ such that, if $i-1 \in D_{i}$, $|D_{i}| = 1$, and $i-2 \notin D_{i+1}, D_{i+2}, \ldots, D_{k}$.  In particular, $(\lambda, A_{1} A_{2} \ldots A_{k}) \in Fix(\widetilde{\psi_{n, d}})$ if and only if $f(\lambda, A_{1} A_{2} \ldots A_{k}) \in Fix(\psi_{n, d})$.
\end{claim}

\begin{proof}
The first part of the claim should be clear from the construction of $f$.  For the second part of the claim, consider the operation of $f$ on a particular $(\lambda, A_{1} A_{2} \ldots A_{k}) \in \widetilde{Y_{n, d}}$: $d$ dots are placed at the tops of their columns, and $X$'s are placed in the columns directed below them and in the rows immediately below them and to their right, with a staircase subdiagram of side length $k$ formed by the remaining free squares.  There will be $\lambda_{j+2}$ rows of $X$'s in between the $j+1$-st and $j+2$-nd rows of this subdiagram, since there will be one row of $X$'s for each dot at the top of its column, and the dots at the top of their columns are placed in blocks of lengths $\lambda_{1}$ (below the first row of the subdiagram), $\lambda_{2}$ (between the first and second rows of the subdiagram) and so on.  

Therefore $j'+1$ precedes $j'$ in $D_{1} D_{2} \ldots D_{k+d} = f(\lambda, A_{1} A_{2} \ldots A_{k})$ for some $j'$ if and only if, for some $j$, $j+1$ precedes $j$ in $A_{1} A_{2} \ldots A_{k}$ and $\lambda_{j+2} = 0$.  Also, $|D_{i}|=1$ for all $1 \leq i \leq k+d$ if and only if $|A_{i}|=1$ for all $1 \leq i \leq k$.  Therefore, $D_{1} D_{2} \ldots D_{k+d} \in Fix(\psi_{n, d})$ if and only if $|A_{i}|=1$ for all $1 \leq i \leq k$ and, for all $p < q$ with $A_{p} = \left\{j+1 \right\}$, $A_{q} = \left\{j \right\}$, $\lambda_{j+2} \neq 0$, or if and only if $f(\lambda, A_{1} A_{2} \ldots A_{k}) \in Fix (\widetilde{\psi_{n, d}})$.  
\end{proof}

\begin{figure}
\includegraphics[width = 10 cm]{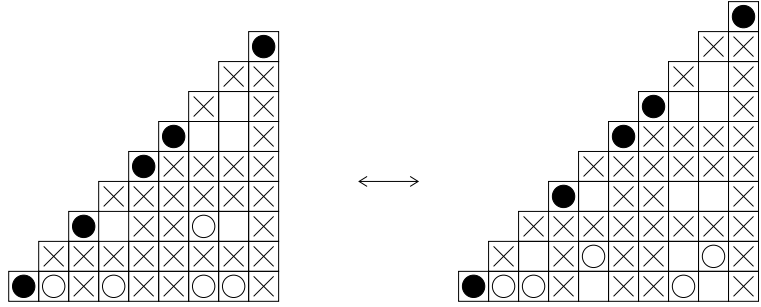}
\caption{$\widetilde{\psi_{10, 5}}: \widetilde{Y_{10, 5}} \rightarrow \widetilde{Y_{10, 5}}$.}
\label{fishburnexam6}
\end{figure}

\begin{figure}
\includegraphics[width = 7 cm]{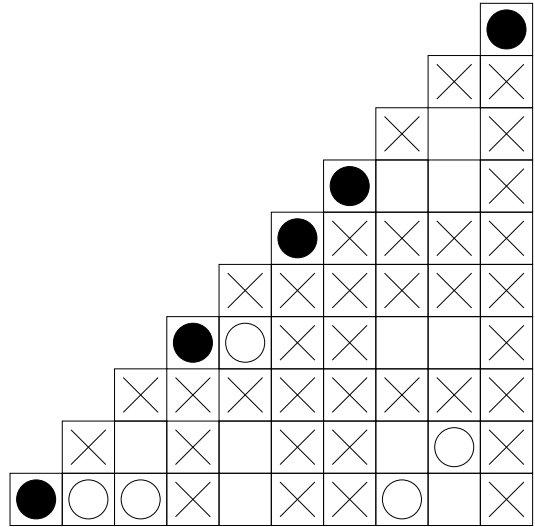}
\caption{An element of $Fix(\widetilde{\psi_{10, 5}})$}
\label{fishburnexam7}
\end{figure}

We illustrate this in Figure \ref{fishburnexam6} and Figure \ref{fishburnexam7}.  The Fishburn diagram on the left-hand-side of Figure \ref{fishburnexam6} corresponds to $f(\lambda, A_{1}A_{2}A_{3}A_{4})$, where $\lambda = (1, 1, 2, 0, 1)$, $A_{1} = A_{2} = A_{4} = \left\{0 \right\}$, and $A_{3} = \left\{0, 1 \right\}$.  Therefore $(\lambda, A_{1} A_{2} A_{3} A_{4}) \in \widetilde{Y_{10, 5}}$.  However,  $(\lambda, A_{1} A_{2} A_{3} A_{4}) \notin Fix(\widetilde{\psi_{10, 5}})$, since $|A_{3}| > 1$, and $\widetilde{\psi_{10, 5}}(\lambda, A_{1} A_{2} A_{3} A_{4}) = (\mu, B_{1} B_{2} B_{3} B_{4} B_{5})$, with $\mu = (1, 0, 1, 2, 0, 1)$, $B_{1} = B_{2} = B_{4} = \left\{0 \right\}, B_{3} = B_{5} = \left\{1 \right\}$.  Note that $1 \in B_{3}, 0 \in B_{4}$, and $0$ is the second element of $\mu$.

By contrast, the Fishburn diagram in Figure \ref{fishburnexam7} corresponds to $f(\mu, E_{1}E_{2}E_{3}E_{4}E_{5})$, where $\mu$ is again equal to $(1, 0, 1, 2, 0, 1)$, $E_{1} = E_{2} = E_{4} = \left\{0 \right\}, E_{3} = \left\{2 \right\}$, $E_{5} = \left\{1 \right\}$.  Although $2 \in E_{3}$ and $1 \in E_{5}$, the third element of $\mu$ is not equal to zero, and so $(\mu, E_{1} E_{2} E_{3} E_{4} E_{5}) \in Fix(\widetilde{\psi_{10, 5}})$.  As we can see, $f(\mu, E_{1}E_{2}E_{3}E_{4}E_{5}) \in Fix(\psi_{10, 5})$.  

Note that $\widetilde{\psi_{n, d}} \neq f^{-1} \circ \psi_{n, d} \circ f$, i.e., $\widetilde{\psi_{n, d}}$ is not simply an alternative description of $\psi_{n, d}$.  In fact, $f^{-1} \circ \psi_{n, d} \circ f$ is not even well-defined, since applying $\psi_{10, 5}$ to the Fishburn diagram on the left-hand-side of Figure \ref{fishburnexam6} will give a Fishburn diagram outside of the image of $f$.  
\section{Further Research Directions}
Our initial approach to proving that matchings with no $2$-nestings were enumerated by the Fishburn numbers was to try to find a bijection between such matchings and matchings with no left- or right-nestings.  Our second approach was to find a bijection between matchings with no $2$-nestings and ascent sequences--in particular, a bijection between matchings $X$ on $[2n]$ on with no $2$-nestings 
such that the inversion table $\phi(X)$ used $k$ distinct integers and ascent sequences of length $n$ with $k-1$ non-ascents. (Inductively, this would make sense: An inversion table of length $n$, with no decreasing subsequence of consecutive integers, that used $k$ distinct integers could have $n+2-k$ integers added to the end to obtain an inversion table of length $n+1$ with no decreasing subsequence of consecutive integers.  An ascent sequence of length $n$ with $k-1$ non-ascents would have $n+1-k$ ascents, so it could have $n+2-k$ integers added to the end to obtain an ascent sequence).  Presumably, bijections between the inversion tables discussed in this paper and other interpretations of the Fishburn numbers could be found by comparing the involutions used to prove the various interpretations have the generating function (\ref{fishburndef}) of the Fishburn numbers, as discussed in Stanley \cite{stanbook}.  It would be interesting to see if any more direct or elegant bijections are possible.

A second interesting question is whether or not it is possible to define bounce, area, or dinv statistics on the various sets enumerated by the Fishburn numbers in such a way as to generalize the bounce, area, and dinv statistics defined on the Catalan sets (see Haglund \cite{hagbook}).  In particular, it would be extremely interesting if any such statistics could be defined so bounce and area or area and dinv gave a symmetric bivariate generating function in general, as the combinatorics of the symmetry of the $q, t$-Catalan polynomial is famously poorly-understood.  Inspired by our proof of the equality of (\ref{refineone}) and (\ref{refinetwo}), we were able to define statistics similar to area and dinv on the inversion tables discussed in this paper as well as slightly different statistics on ascent sequences.  In both cases we obtained symmetric generating functions for all $n \leq 5$, but not for $n = 6$.  (We thank Jason Bandlow for his assistance with checking the $n=6$ cases by computer).  The polynomial obtained from ascent sequences was closer to being symmetric than the polynomial obtained from inversion tables with no decreasing subsequence of consecutive integers; in both cases, $f(q, t) - f(t, q)$ was only a few terms long.  It is possible that working with another interpretation of the Fishburn numbers would inspire statistics that were, in fact, symmetric.

\section{Appendix}
\begin{figure}
\includegraphics[width = 10 cm]{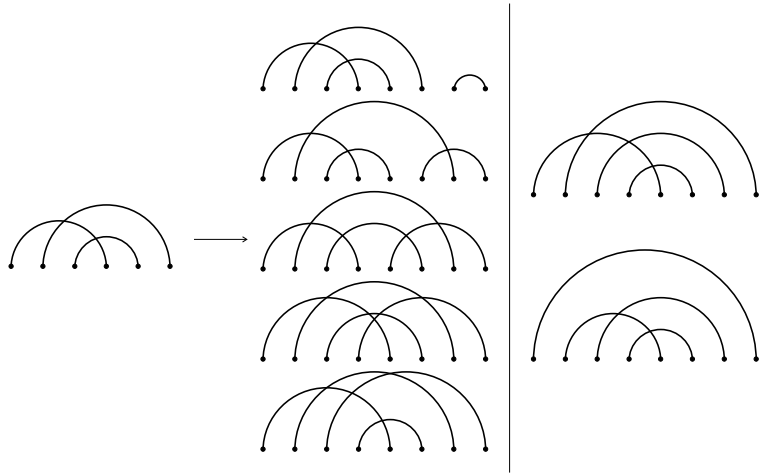}
\caption{Matchings on $[6]$ and $[8]$}
\label{eulerexam1}
\end{figure}

In this Appendix, we will prove the following conjecture of Claesson and Linusson's:
\begin{claim} The distribution of left-nestings over the set of all matchings on $[2n]$ is
given by the ``Second-order Eulerian triangle'', entry $A008517$ in OEIS \cite{OEIS}.
\end{claim}
\begin{proof}
According to the OEIS, the Second-order Eulerian triangle $T(n, k)$ is defined by the following recurrence relation: $T(n,k) = 0$ if $n<k$, $T(1,1)=1$, $T(n,-1)=0$, $T(n,k)=kT(n-1,k)+(2n-k)T(n-1,k-1).$

Assume $X$ is a matching on $[2n-2]$ with precisely $j$ left-nestings.  If we add an arc corresponding to a pair $(a, 2n)$ to achieve a matching $X'$ on $[2n]$, there are three possibilities:
\begin{description}
\item[Case 1] If the left-endpoint of this new arc is placed immediately before a right-endpoint of $X$, or immediately before its own right-endpoint, this will not change the number of left-nestings, so $X'$ will have $j$ left-nestings.  There are $n$ right-endpoints (including $2n$ itself), so there are $n$ ways to do this.
\item[Case 2] If the left-endpoint of this new arc is placed immediately before a left-endpoint of an arc that is \textit{not} the inner arc of a left-nesting of $X$, this will contribute one new left-nesting to $X'$, without eliminating an existing left-nesting of $X$.  Then $X'$ will have $j+1$ left-nestings, and there are $n-1-j$ ways to do this.
\item[Case 3] If the left-endpoint of this new arc is placed immediately before a left-endpoint of an arc that \textit{is} the inner arc of a left-nesting of $X$, this will add one new left-nesting to $X'$, but also eliminate the original left-nesting from $X$.  The net effect is therefore to contribute no new left-nestings to $X'$, which will therefore have $j$ left-nestings.  There are $j$ ways to do this.
\end{description}
Therefore, of the $2n-1$ matchings on $[2n]$ formed by adding a rightmost arc to $X$, $n+j$ will have $j$ left-nestings, and $n-1-j$ will have $j+1$ left-nestings. 
To put it another way, if $L(n, j)$ is defined to be the number of matchings on $[2n]$ with $j$ left-nestings, we now have that:
\begin{eqnarray}
L(n, j) = (n+j)L(n-1, j) + (n-j)L(n-1, j-1),
\end{eqnarray}

or, replacing $j$ with $n-k$,
\begin{eqnarray}
L(n, n-k) = (2n-k)L(n-1, n-1-(k-1)) + kL(n-1, n-1-k).
\end{eqnarray}
Since $L(1, 1-1) = 1$, $L(n, n-(-1)) = 0$, and $L(n, n-k) = 0$ if $n < k$, we see inductively that $L(n, n-k) = T(n, k).$  
\end{proof}
We illustrate the three cases in the above inductive proof in Figure \ref{eulerexam1}, beginning with the matching $(1, 4) (2, 6) (3, 5)$ on $[6]$, which has $1$ left-nesting.  Adding the arcs corresponding to the pairs $(7, 8), (6, 8), (5, 8), (4, 8)$ and re-labeling falls into Case $1$ and results in matchings on $[8]$ with $1$ left-nesting.  Adding the arc corresponding to the pair $(3, 8)$ and re-labelling falls into case $3$ and also results in a matching on $[8]$ with $1$ left-nesting.  Adding the arcs corresponding to the pairs $(2, 8)$ and $(1, 8)$ and re-labelling falls into Case $3$ and results in a matching on $[8]$ with $2$ left-nestings.  
Note that the above proof can be seen as a generalization of the inductive proof of Claim (\ref{inductclaim}). 

\section{Acknowledgements}
The author would like to thank Anders Claesson, Svante Linusson, and Jeffrey Remmel for their encouragement and advice throughout this research.  
\bibliography{mybib}{}
\bibliographystyle{hplain}
\end{document}